\UseRawInputEncoding
\documentclass[12pt,reqno]{amsart}

\usepackage[english]{babel} 
\usepackage{amsmath,amsfonts,amssymb,bm} 
\usepackage{verbatim}
\usepackage{mathrsfs}

\newtheorem{thm}{Theorem}
\newtheorem{prop}{Proposition}
\newtheorem{cor}{Corollary}
\newtheorem{lem}{Lemma}
\newtheorem{citethm}{Theorem}

\theoremstyle{definition}
\newtheorem{rem}{Remark}
\newtheorem{ex}{Example}

\providecommand{\RR}{\mathbb{R}}

\providecommand{\HH}{\mathbb{H}}

\providecommand{\NN}{\mathbb{N}}

\DeclareMathOperator{\id}{id}

\DeclareMathOperator{\im}{Im}
\DeclareMathOperator{\spa}{span}

\def\bphi{\boldsymbol \phi}

\makeatletter
\@namedef{subjclassname@2020}{%
  \textup{2020} Mathematics Subject Classification}
\makeatother

\begin{document}

\title[Matrix version of Quaternionic Nullstellensatz]{Matrix versions of real and quaternionic nullstellensatz}

\begin{abstract}
Real Nullstellensatz is a classical result from Real Algebraic Geometry.
It has recently been extended to quaternionic polynomials by Alon and 
Paran \cite{qsatz1}. The aim of this paper is to extend their
Quaternionic Nullstellensatz to matrix polynomials. 
We also obtain an improvement of the Real Nullstellensatz for matrix polynomials
from \cite{c1} in the sense that we simplify the definition of a real left ideal.
We use the methods from the proof of the matrix version of Hilbert's Nullstellensatz \cite{c2} and we obtain their extensions to a mildly non-commutative case and to the real case.
\end{abstract}

\author{J. Cimpri\v c}

\date{January 1st 2022}

\thanks{This work was supported by grant P1-0222 from Slovenian Research Agency}



\address{University of Ljubljana, Faculty of Mathematics and Physics, Department of Mathematics, Jadranska 21, 1000 Ljubljana, Slovenia}

\email{jaka.cimpric@fmf.uni-lj.si}

\maketitle

\thispagestyle{empty}

\section{Introduction}

\subsection{Real Nullstellensatz}
The aim of this paper is to discuss several non-commutative generalizations
of the following result:

\begin{citethm}[Real Nullstellensatz, \cite{dub,ris,efr}]
\label{thma}
For any polynomials $p_1,\ldots,p_m,q \in \RR[x_1,\ldots,x_d]$ the following are equivalent:
\begin{enumerate}
\item For every point $a \in \RR^d$ such that $p_1(a)=\ldots=p_m(a)=0$, we have $q(a)=0$.
\item $q$ belongs to the smallest real ideal of $\RR[x_1,\ldots,x_d]$ that contains $p_1,\ldots,p_m$.
\end{enumerate}
\end{citethm}

Recall that an ideal $I$ of a commutative ring $R$ is \textit{real} if for every
$a_1,\ldots,a_k \in R$ such that $a_1^2+\ldots+a_k^2 \in I$
we have $a_1,\ldots,a_k \in I$. The smallest real ideal that contains
a given ideal $J$ is called the \textit{real radical} of $J$ and it is
usually denoted by $\sqrt[\text{rr}]{J}$.

\subsection{Quaternionic polynomials}
In \cite{qsatz1}, Theorem \ref{thma} was extended to quaternionic polynomials.
Their results are summarized in subsections \ref{sub14} and \ref{sub15}.
They work with two types of quaternionic polynomials. The first one is
the ring of polynomial functions from $\RR^d$ to $\HH$, i.e.
$$\HH_c[x_1,\ldots,x_d] := \HH \otimes_\RR \RR[x_1,\ldots,x_d]$$
whose variables commute with each other and with the elements of $\HH$. 
Its involution is defined by $x_i^\ast=x_i$ for all $i$ and $c^\ast=\bar{c}$
for all $c \in \HH$. The second one is the ring of all polynomial functions from
$\HH^d$ to $\HH$;
$$\HH[x_1,\ldots,x_d] := \im \left(\HH\langle x_1,\ldots,x_d \rangle
\to \text{Functions}(\HH^d,\HH) \right)$$ 
whose variables do not commute with each other or the elements of $\HH$.
Note that the ring $\HH[x_1,\ldots,x_d]$ is invariant under the pointwise conjugation on 
$\text{Functions}(\HH^d,\HH)$ since $\bar{f}=-\frac12(f+ifi+jfj+kfk)$.
This defines a natural involution on $\HH[x_1,\ldots,x_d]$.

The aim of this paper is to extend the results of \cite{qsatz1} to matrices over quaternionic polynomials.
Our new results are summarized in subsections \ref{sub16}, \ref{sub17}, \ref{sub18} and \ref{sub19}. Reference \cite{qsatz2} is a different generalization of \cite{qsatz1} that will not be discussed here. References \cite{qsatz3} and \cite{ones} are earlier developments. 

\subsection{Real left ideals}
\label{sub13}
Let $A$ be a ring with involution. We say that a left ideal $I$ of $A$ is \textit{real}
if for every $a_1,\ldots,a_k \in A$ such that $a_1^\ast a_1+\ldots+a_k^\ast a_k \in I$
we have $a_1,\ldots,a_k \in I$.
 
\begin{rem}
\label{defrem1}
If the involution of $A$ is trivial,  we get the usual definition 
of a real ideal in a commutative ring. 
\end{rem}

\begin{rem}
\label{defrem2}
If $a^\ast a=a a^\ast$ for every $a \in A$ (e.g. for quaternionic polynomials) 
then every real left ideal is $\ast$-invariant and therefore two-sided.
\end{rem}

\begin{rem}
\label{defrem4}
An intersection of a family of real left ideals is clearly a real left ideal.
For every subset of $A$ there therefore exists the smallest real left
ideal that contains it (i.e. its real radical.)
\end{rem}

\begin{rem}
\label{defrem5}
Every real left ideal $I$ of $A$ is semiprime. Namely, if $xAx \subseteq I$ for some 
$x \in A$ then $x^\ast xx^\ast x \in I$. If $I$ is real, it follows that $x^\ast x \in I$ 
and so $x \in I$.
\end{rem}

\begin{rem}
\label{defrem3}
Our definition of a real left ideal is different from the definition in \cite{chmn} 
or \cite{c1}. In subsection \ref{sub19} we will show that for the rings discussed
in this paper both definitions are equivalent. We do not know if they are equivalent 
in general.
\end{rem}

\subsection{First Quaternionic Nullstellensatz}
\label{sub14}
Theorem \ref{thmb} is a generalization of Theorem \ref{thma} 
to $\HH_c[x_1,\ldots,x_d]$. Its proof is hidden in the proof of \cite[Theorem 9]{qsatz1}
but we summarize it here for the sake of completeness.

\begin{citethm}
\label{thmb}
For every $p_1,\ldots,p_m,q \in \HH_c[x_1,\ldots,x_d]$ the following are equivalent:
\begin{enumerate}
\item For every point $a \in \RR^d$ such that $p_1(a)=\ldots=p_m(a)=0$, we have $q(a)=0$.
\item $q$ belongs to the smallest real ideal of $\HH_c[x_1,\ldots,x_d]$ that contains $p_1,\ldots,p_m$.
\end{enumerate}
\end{citethm}

\begin{proof}
 Suppose that claim (1) of Theorem \ref{thmb} is true.	
Write $p_i=p_{i,1} \mathbf{1}+p_{i,2} \mathbf{i}+p_{i,3} \mathbf{j}+p_{i,4} \mathbf{k}$, $q=q_1 \mathbf{1}+q_2 \mathbf{i}+q_3 \mathbf{j}+q_4 \mathbf{k}$
and let $\mathfrak{a}$ be the smallest real ideal of $R=\RR[x_1,\ldots,x_d]$ that contains all $p_{i,j}$.
By Theorem \ref{thma}, it follows that $q_k \in \mathfrak{a}$ for all $k$ and so $q \in \HH \cdot \mathfrak{a}$. 
Finally, note that $\HH \cdot \mathfrak{a}$ is the smallest real two-sided ideal of 
$\HH_c[x_1,\ldots,x_d]=\HH \otimes R$ that contains $p_1,\ldots,p_m$.
Namely, by \cite[Lemma 17.4]{gwbook}, the mappings $\iota \mapsto \HH \cdot \iota$ 
and $I \mapsto I \cap R$ give a 1-1 correspondence between the ideals of $R$ 
and the two-sided ideals of $\HH \otimes R$. 
Clearly this correspondence preserves reality; see Remark \ref{preserve}.
Therefore, we have claim (2) of Theorem \ref{thmb}.
\end{proof}

\subsection{Second Quaternionic Nullstellensatz}
\label{sub15}
Theorem \ref{thmc} is a generalization of Theorem \ref{thma} to $\HH[x_1,\ldots,x_d]$.
It rephrases \cite[Theorem 9]{qsatz1} in our terminology.

\begin{citethm}
\label{thmc}
For every $p_1,\ldots,p_m,q \in \HH[ x_1,\ldots,x_d ]$ the following are equivalent
\begin{enumerate}
\item For every $a \in \HH^d$ such that $p_1(a)=\ldots=p_m(a)=0$ we have $q(a)=0$.
\item $q$ belongs to the smallest real ideal of $\HH[ x_1,\ldots,x_d ]$ containing $p_1,\ldots,p_m$.
\end{enumerate}
\end{citethm}

\begin{proof}
By \cite[Theorem 5]{qsatz1}, the mapping
\begin{eqnarray}
\label{mainiso}
\phi \colon \HH[x_1,\ldots,x_d] & \to & 
\HH_c[y_{1,1},\ldots,y_{1,4},\ldots, y_{d,1},\ldots,y_{d,4}]\\ \notag
x_i & \mapsto & 
y_{i,1}\mathbf{1}+y_{i,2}\mathbf{i}+y_{i,3} \mathbf{j}+y_{i,4} \mathbf{k}
\end{eqnarray}
is well-defined and it is an isomorphism of rings with involution. 
We identify quaternionic points with real points by the natural mapping 
$\rho \colon \HH^d \to \RR^{4d}$ and we observe that 
$f(a)=\phi(f)(\rho(a))$ for every $f \in \HH[x_1,\ldots,x_d]$
and every $a \in \HH^d$.

Suppose that claim (1) of Theorem \ref{thmc} is true.
It follows that for every $a \in \HH^d$ which satisfies 
$\phi(p_i)(\rho(a))=0$ for every $i$ we have $\phi(q)(\rho(a))=0$.
By Theorem \ref{thmb}, it follows that $\phi(q)$ belongs to the
smallest real ideal of $\HH_c[y_{1,1},\ldots,y_{d,4}]$
which contains $\phi(p_1),\ldots,\phi(p_m)$. Since $\phi$ is
a $\ast$-isomorphism, it follows that claim (2) of Theorem \ref{thmc} is true.
\end{proof}

\subsection{Matricial Real Nullstellensatz}
\label{sub16}
The aim of this paper is to prove matrix versions of Theorems \ref{thma},
\ref{thmb} and \ref{thmc}. For every ring with involution $A$,
we will assume that the ring $M_n(A)$ has the conjugated transpose 
involution, i.e. $[a_{i,j}]^\ast=[a_{j,i}^\ast]$.

Theorem \ref{thmd} is a matrix version of Theorem \ref{thma}.
Its proof uses the main results from \cite{c2}
which are summarized in Section \ref{sec2}. 
Real versions of these results are developed in Section \ref{sec3}
and they are used to prove a vector version of Theorem \ref{thma}; see Corollary \ref{realthm2}.
Finally, the proof of Theorem \ref{thmd} is given in Section \ref{sec4}.

\begin{citethm}[Matricial Real Nullstelensatz]
\label{thmd}
For any real matrix polynomials $P_1,\ldots,P_m,Q \in M_n(\RR[x_1,\ldots,x_d])$ the following are equivalent:
\begin{enumerate}
\item For every $a \in \RR^d$ and every $\mathbf{v} \in \RR^n$ such that 
$P_1(a)\mathbf{v}=\ldots=P_m(a)\mathbf{v}=0$, we have $Q(a)\mathbf{v}=0$.
\item $Q$ belongs to the smallest real left ideal of $M_n(\RR[x_1,\ldots,x_d])$ that contains $P_1,\ldots,P_m$.
\end{enumerate}
\end{citethm}

Theorem \ref{thmd} implies the main result of \cite{c1}
and  the equivalence of our definition of a real left ideal
to the definition in \cite{c1}; see Corollary \ref{eqdef}.

\subsection{First Matricial Quaternionic Nullstellensatz}
\label{sub17}
Theorem \ref{thme} is a matrix version of Theorem \ref{thmb}.
Its proof is spread out over sections \ref{sec5}, \ref{sec6} and \ref{sec7}.
We extend first the results of sections \ref{sec2} and \ref{sec3}
to a class of non-commutative rings that includes quaternionic polynomials.
Then we  prove a vector version of Theorem \ref{thmb}; see Corollary \ref{qthm2}.
Finally, Theorem \ref{thme} is deduced from Corollary \ref{qthm2}.

\begin{citethm}
\label{thme}
For any quaternionic matrix polynomials  $P_1,\ldots,P_m,Q \in M_n(\HH_c[x_1,\ldots,x_d])$ the following are equivalent:
\begin{enumerate}
\item For every $a \in \RR^d$ and every $\mathbf{v} \in \HH^n$ such that $P_1(a)\mathbf{v}=\ldots=P_m(a)\mathbf{v}=0$, we have $Q(a)\mathbf{v}=0$.
\item $Q$ belongs to the smallest real left ideal of $M_n(\HH_c[x_1,\ldots,x_d])$ that contains $P_1,\ldots,P_m$.
\end{enumerate}
\end{citethm}

\subsection{Second Matricial Quaternionic Nullstellensatz}
\label{sub18}
Theorem \ref{thmf} is a matrix version of Theorem \ref{thmc}.
We deduce it from Theorem \ref{thme} by means of isomorphism \eqref{mainiso}
(analogously to the proof of Theorem \ref{thmc}.)

\begin{citethm}
\label{thmf}
For any quaternionic matrix polynomials  $P_1,\ldots,P_m,Q \in M_n(\HH [x_1,\ldots,x_d ])$ the following are equivalent:
\begin{enumerate}
\item For every $a \in \HH^d$ and $\mathbf{v} \in \HH^n$ such that $P_1(a)\mathbf{v}=\ldots=P_m(a)\mathbf{v}=\mathbf{0}$ we have $Q(a)\mathbf{v}=\mathbf{0}$.
\item $Q$ belongs to the smallest real left ideal of $M_n(\HH [x_1,\ldots,x_d ])$ containing $P_1,\ldots,P_m$.
\end{enumerate}
\end{citethm}

\begin{proof}
The $\ast$-isomorphism $\phi$ from \eqref{mainiso} induces an $\ast$-isomorphism 
$$\phi_n \colon M_n(\HH[x_1,\ldots,x_d]) \to M_n(\HH_c [y_{1,1},\ldots,y_{1,4},\ldots,y_{d,1},\ldots,y_{d,4}])$$
which satisfies $F(a)=\phi_n(F)(\rho(a))$
for every $F \in M_n(\HH[x_1,\ldots,x_d])$ and every $a \in \HH^d$.
(Recall that $\rho \colon \HH^d \to \RR^{4d}$ is the natural map.)

Claim (1) of Theorem \ref{thmf} implies that for every 
$a \in \HH^d$ and $\mathbf{v} \in \HH^n$ such that 
$\phi_n(P_1)(\rho(a))\mathbf{v}=\ldots=\phi_n(P_m)(\rho(a))\mathbf{v}=0$ 
we have $\phi_n(Q)(\rho(a))\mathbf{v}=0$.
By Theorem \ref{thme}, it follows that $\phi_n(Q)$ belongs to
the smallest real left ideal of $M_n(\HH_c [y_{1,1},\ldots,y_{1,4},\ldots,y_{d,1},\ldots,y_{d,4}])$ that contains $\phi_n(P_1),\ldots,\phi_n(P_m)$.
Since $\phi_n$ is a $\ast$-isomorphism, this implies claim (2) of Theorem \ref{thmf}.
\end{proof}

\subsection{An equivalent definition of a real left ideal}
\label{sub19}
Let us compare our definition of a real left ideal with the definition in \cite{chmn} or
\cite{c1}. We say that a left ideal $I$ of a ring with involution $A$ is 
\textit{strongly real} (i.e. real in the sense of \cite{chmn}, \cite{c1})
if for every $a_1,\ldots,a_k \in A$ such that 
$a_1^\ast a_1+\ldots+a_k^\ast a_k \in I+I^\ast$ we have $a_1,\ldots,a_k \in I$.

\begin{cor}
\label{eqdef}
Let $A$ be either $M_n(\RR[x_1,\ldots,x_d])$ or $M_n(\HH_c[x_1,\ldots,x_d])$
or $M_n(\HH[x_1,\ldots,x_n])$ with the usual involutions. A left ideal of $A$
is real iff it is strongly real.
\end{cor}

\begin{proof}
Let $A=M_n(\HH[x_1,\ldots,x_n])$. For other cases the proof is similar.

We will prove first that for every subset $S$ of $\HH^d \times \HH^n$ the left ideal
$$J(S)=\{F \in A \mid F(a)\mathbf{v}=\mathbf{0} \text{ for all } 
(a,\mathbf{v}) \in S\}$$
is strongly real. Namely, if $\sum_{i=1}^m F_i^\ast F_i$ belongs to $J(S)+J(S)^\ast$
for some $F_1,\ldots,F_m \in A$ then there exist $G,H \in J(S)$ such that 
$\sum_{i=1}^m F_i^\ast F_i =G+H^\ast$. For every $(a,\mathbf{v}) \in S$
we have $G(a)\mathbf{v}=\mathbf{0}$ and $H(a)\mathbf{v}=\mathbf{0}$
which implies that $\sum_{i=1}^m \mathbf{v}^\ast F_i(a)^\ast F_i(a) \mathbf{v}=0$.
The latter implies that $F_i(a)\mathbf{v}=\mathbf{0}$ for all $i$, 
i.e. $F_i \in J(S)$ for all $i$.

On the other hand, by Theorem \ref{thmf} and the well-known fact that $A$ is left Noetherian,
every real left ideal $I$ of $A$
is equal to $J(V(I))$ where 
$$V(I)=\{(a,\mathbf{v}) \in \HH^d \times \HH^n \mid F(a)\mathbf{v}=\mathbf{0}
\text{ for all } F \in I\}$$
Therefore every real left ideal of $A$ is strongly real. The opposite direction is clear.
\end{proof}

\section{Prime and semiprime submodules}
\label{sec2}

In this section we summarize the main results 
from \cite{c2} which will be needed for the proof of Theorem \ref{thmd}.
We also make preparations for the proofs of mildly non-commutative versions of these results 
in Sections \ref{sec5} and \ref{sec6} (which are needed for the proof of Theorem \ref{thme}.)

Let $A$ be an associative ring with $1$ and $n$ a positive integer.
We will write $S$ for the full matrix ring $M_n(A)$ and $M$ for the left $A$-module $A^n$.
There is a natural 1-1 correspondence between the left ideals of $S$ 
and the submodules of $M$. 
(The left ideal of $S$ that corresponds to a submodule $N$ of $M$ consists of all matrices
that have all rows in $N$.)
 There is also a natural 1-1 correspondence 
between the two-sided ideals of $S$ and the two-sided ideals of $A$. 
(The two-sided ideal of $S$ that corresponds to a two-sided ideal $J$ of $A$
consists of all matrices that have all entries in $J$.)
For every submodule $N$ of $M$ the set
$$(N:M):=\{a \in A \mid aM \subseteq N\}$$ is a two-sided ideal of $A$.
The two-sided ideal of $S$ that corresponds to
$(N:M)$ is equal to the largest two-sided ideal  contained 
in the left ideal of $S$ that corresponds to $N$.

A submodule $N$ of $M$ is \textit{semiprime} 
iff for every $\mathbf{f}=(f_1,\ldots,f_n) \in M$ 
such that $f_i A \mathbf{f} \subseteq N$ 
for all $i$, we have that $\mathbf{f} \in N$.
If $N$ is a semiprime submodule of $M$ then 
$(N:M)$ is clearly a semiprime two-sided ideal of $A$.
Moreover, by \cite[Proposition 3]{c2}, the left ideal of $S$
that corresponds to $N$ is also semiprime. Recall that
a left ideal $I$ of $S$ is semiprime if for every $x \in S$ 
such that $xSx \subseteq I$ we have $x \in I$. 

A submodule $N$ of $M$ is \textit{prime} if for every $a \in A$ and 
every $\mathbf{m} \in M$ such that $a A \mathbf{m} \subseteq N$ 
we have that either $a \in (N:M)$ or $\mathbf{m} \in N$. 
If $N$ is a prime submodule of $M$ then $(N:M)$ is a prime
two-sided ideal of $A$ and, by \cite[Proposition 3]{c2},
the left ideal of $S$ that corresponds to $N$ is also prime.
Recall that a left ideal $I$ of $S$ is prime if for every $x,y \in S$
such that $xSy \subseteq I$ we have either $xS \subseteq I$ or $y \in I$.

Let $\mathfrak{p}$ be a two-sided prime ideal of $A$. A prime submodule $N$ of $M$ such that $(N : M)=\mathfrak{p}$ is called $\mathfrak{p}$-\textit{prime}. 

\begin{rem}
	\label{altpp}
If $\mathfrak{p}=A$ then a submodule $N$ of $M$ is $\mathfrak{p}$-prime iff $N=M$.
If $\mathfrak{p} \ne A$ then a submodule $N$ of $M$ is $\mathfrak{p}$-prime iff
$N \ne M$ and $\mathfrak{p}^n \subseteq N$ and 	
for every $a \in A$ and $\mathbf{m} \in M$ such that $a A\mathbf{m} \subseteq  N$
we have either $a \in \mathfrak{p}$ or $\mathbf{m} \in N$. Namely, the second condition implies
that $\mathfrak{p} \subseteq (N : M)$ while the first and the third condition imply that
$\mathfrak{p} \supseteq (N : M)$.
\end{rem}
	
\begin{ex} 
	\label{smallest}
The smallest $\mathfrak{p}$-prime submodule of $M$ is $\mathfrak{p}^n$.
\end{ex}

\begin{ex}
\label{cpu}
For every $\mathbf{u} \in M\setminus \mathfrak{p}^n$, the set 
$$C_{\mathfrak{p},\mathbf{u}}:=\{ \mathbf{m} \in M \mid 
\langle \mathbf{m},\mathbf{u} \rangle \in \mathfrak{p}\}$$ 
(where $\langle \mathbf{m},\mathbf{u} \rangle = \sum_{i=1}^n m_i u_i$)
is a $\mathfrak{p}$-prime submodule of $M$. 
Namely, pick any $a \in A$ and $\mathbf{m} \in M$ such that 
$a A\mathbf{m} \subseteq C_{\mathfrak{p},\mathbf{u}}$.
It follows that for every $x \in A$, $a x \langle \mathbf{m},\mathbf{u} \rangle 
=\langle a x \mathbf{m},\mathbf{u} \rangle \in \mathfrak{p}$
which implies that either $a \in \mathfrak{p}$ or $\mathbf{m} \in C_{\mathfrak{p},\mathbf{u}}$.
Now use Remark \ref{altpp}.
\end{ex}

\begin{ex}
\label{knp}
An intersection of a family of $\mathfrak{p}$-prime submodules is again a 
$\mathfrak{p}$-prime submodule by Remark \ref{altpp}. If a submodule $N$ of $M$ is 
contained in some $\mathfrak{p}$-prime submodule we write 
$K(N,\mathfrak{p})$ for the smallest $\mathfrak{p}$-prime submodule 
that contains it. Otherwise we write $K(N,\mathfrak{p})=M$.

If $A$ is commutative and $\mathfrak{p} \ne A$ then 
$K(N,\mathfrak{p}) \ne M$ iff $\mathfrak{p} \supseteq (N : M)$ by claim (4) of Proposition \ref{class}.
This is also true for some non-commutative $A$; see Theorem \ref{ncintthm2}. 
For general $A$, however, we can only show that $K(N,\mathfrak{p}) \ne M$ 
implies $\mathfrak{p} \supseteq (N : M)$. Namely, if $K(N,\mathfrak{p}) \ne M$
and $(N : M) \not\subseteq \mathfrak{p}$ then $K(N,\mathfrak{p})$ is $\mathfrak{p}$-prime
and we can  pick $x \in (N : M)$ such that $x \not\in \mathfrak{p}$.
Since $x \not\in \mathfrak{p}$ and $x A \mathbf{m} \subseteq N \subseteq K(N,\mathfrak{p})$ 
for every $\mathbf{m} \in M$, it follows that $\mathbf{m} \in K(N,\mathfrak{p})$ for every $\mathbf{m} \in M$; a contradiction. 
\end{ex}

\begin{prop}
\label{class}
Let $R$ be a commutative ring and let $n\in \NN$.
Pick a prime ideal $\mathfrak{p} \ne R$ and write  $k$ for the field of fractions of $R/\mathfrak{p}$.
\begin{enumerate}
\item The natural mapping $\pi \colon R^n \to k^n$ induces a 1-1 correspondence
between $\mathfrak{p}$-prime submodules of $R^n$ and proper subspaces of $k^n$.
For convenience we will also say that $R^n$ corresponds to $k^n$.
\item A submodule of the form $C_{\mathfrak{p},\mathbf{u}}$ corresponds to the
maximal proper subspace ``orthogonal'' to $\pi(\mathbf{u})$. Therefore, submodules of the form $C_{\mathfrak{p},\mathbf{u}}$ coincide with maximal $\mathfrak{p}$-prime submodules of $R^n$.
\item For every submodule $N$ of $R^n$, the submodule $K(N,\mathfrak{p})$ corresponds to the subspace $\spa \phi(N)$. It follows that
\begin{equation}
\label{knpf}
K(N,\mathfrak{p})=\{\mathbf{m} \in R^n \mid \exists c \in R \setminus \mathfrak{p} \colon c  \mathbf{m} \in N+\mathfrak{p}^n\}.
\end{equation}
\item For every submodule $N \ne R^n$ we have $K(N,\mathfrak{p}) \ne R^n$ iff
$\mathfrak{p}$ contains $(N : R^n)$. In this case, we have the following formula
\begin{equation}
\label{kcf}
K(N,\mathfrak{p})=\bigcap_{N \subseteq C_{\mathfrak{p},\mathbf{u}}} 
C_{\mathfrak{p},\mathbf{u}}
\end{equation}
\end{enumerate}
\end{prop}

\begin{proof}
For claims (1) and (3) see \cite[Proposition 1]{c2}.
For claim (2) and formula \eqref{kcf} see \cite [Corollary 1]{c2}.

To prove the first part of claim (4) suppose that $K(N,\mathfrak{p})=M$. 
Let $\mathbf{e}_1,\ldots,\mathbf{e}_n$ be the standard basis of $M$. 
By formula \eqref{knpf} there exist 
$c_1,\ldots,c_n \in R \setminus \mathfrak{p}$ such that
$c_i \mathbf{e}_i \in N+\mathfrak{p}^n$ for all $i$. 
Write $c=c_1 \cdots c_n \in R \setminus \mathfrak{p}$
and note that for every $\mathbf{m}=\sum_i m_i \mathbf{e}_i \in M$ 
we have $c \mathbf{m}=\sum_i m_i c \mathbf{e}_i \in N+\mathfrak{p}^n$. 
If we consider $M$ as a subset of $(R_{\mathfrak{p}})^n$, we get
$M \subseteq c^{-1} (N+\mathfrak{p}^n) \subseteq R_{\mathfrak{p}} N+\mathfrak{p}
(R_{\mathfrak{p}})^n$. It follows that $(R_{\mathfrak{p}})^n=R_{\mathfrak{p}}N+
\mathrm{Jac}(R_{\mathfrak{p}})(R_{\mathfrak{p}})^n$. By Lemma Nakayama, it follows that $(R_{\mathfrak{p}})^n=R_{\mathfrak{p}}N$. By clearing denominators,
we obtain $x \in R \setminus \mathfrak{p}$ such that $x M \subseteq N$.
In other words, $(N : M) \not \subseteq \mathfrak{p}$.
For the converse see Example \ref{knp}.
\end{proof}

\begin{thm}
\label{intthm1}
Let $A$ be an associative ring and let $N$ be a submodule of $M=A^n$. Consider the following claims:
\begin{enumerate}
\item $N$ is semiprime, 
\item $N=\bigcap_{\mathfrak{p}} K(N,\mathfrak{p})$
where $\mathfrak{p}$ runs through all two-sided prime ideals of $A$,
\item $N=\bigcap_{N \subseteq C_{\mathfrak{p},\mathbf{u}}}
C_{\mathfrak{p},\mathbf{u}}$ where $\mathfrak{p}$ runs through all two-sided prime ideals of $A$ and $\mathbf{u}$ runs through all elements of $M \setminus \mathfrak{p}^n$.
\item $N$ is equal to the intersection of all prime submodules that contain it.
\end{enumerate}
We always have $(3) \Rightarrow (4) \Rightarrow (1)$.
If $A$ is commutative, we also have $(1) \Rightarrow (2) \Rightarrow (3)$.
\end{thm}

Later we will show that the implications $(1) \Rightarrow (2) \Rightarrow (3)$ 
also hold for a  class of non-commutative rings that contains the ring $\HH_c[x_1,\ldots,x_d]$;
see Theorems \ref{ncintthm2} and \ref{ncmax}.

\begin{proof}
$(3) \Rightarrow (4)$ follows from Example \ref{cpu}.

$(4) \Rightarrow (1)$ follows from the fact that every prime submodule is semiprime
and the fact that the intersection of a family of semiprime submodules is a semiprime
submodule.

Let us assume now that $A$ is commutative.

$(1) \Rightarrow (2)$ is the same as \cite[Theorem 1]{c2}.

$(2) \Rightarrow (3)$ follows from claim (4) of Proposition \ref{class}.
\end{proof}

\section{Real submodules}
\label{sec3}

The aim of this section is to prove Theorem \ref{realthm1} which is a real version of Theorem \ref{intthm1}.
As a corollary we obtain a vector version of Theorem \ref{thma}
which is the crucial step in the proof of Theorem \ref{thmd}.
Later we will also prove a non-commutative version of Theorem
\ref{realthm1}; see Theorems \ref{ncmax} and \ref{qthm1}.

Let $A$ be an associative ring with involution $\ast$ and let $N$ be a submodule of $A^n$.
We say that $N$ is \textit{real} if for every 
$$\mathbf{f}_i=(f_{i,1},\ldots,f_{i,n}) \in A^n \text{ where } i=1,\ldots,m$$
such that
$$\sum_{i=1}^m f_{i,j}^\ast \mathbf{f}_i \in N \text{ for every } j=1,\ldots,n$$
we have that $$\mathbf{f}_i \in N \text{ for every } i=1,\ldots,m.$$
For $n=1$ this definition coincides with the definition
of a real left ideal of $A$; see subsection \ref{sub13}.

We want to prove the following.

\begin{thm}
	\label{realthm1}
	Let $A$ be an associative ring with involution and let $N$ be a submodule of $M=A^n$. Consider the following statements:
	\begin{enumerate}
		\item $N$ is real,
		\item $N$ is semiprime and the ideal $(N : M)$ is real.
		\item $N=\bigcap_{\mathfrak{p}} K(N,\mathfrak{p})$ where $\mathfrak{p}$ runs through all two-sided real prime ideals of $A$.
		\item $N=\bigcap_{N \subseteq C_{\mathfrak{p},\mathbf{u}}}C_{\mathfrak{p},\mathbf{u}}$ where $\mathfrak{p}$ runs through all two-sided real prime ideals of $A$ and $\mathbf{u}$ runs through all elements of $M \setminus \mathfrak{p}^n$.
		\item $N$ is equal to the intersection of all real prime submodules that contain it.
	\end{enumerate}
We always have $(4) \Rightarrow (5) \Rightarrow (1) \Rightarrow (2)$.
	If $A$ is commutative with trivial involution, we also have $(2) \Rightarrow (3) \Rightarrow (4)$.
\end{thm}

Later we will show that the implications $(2) \Rightarrow (3) \Rightarrow (4)$ 
also hold for a class of non-commutative rings with involution that contains $\HH_c[x_1,\ldots,x_d]$;
see Theorems \ref{ncmax} and \ref{qthm1}.

\begin{proof}
	$(1) \Rightarrow (2)$. Let $N$ be a real submodule of $M:=A^n$.
	
	To prove that $N$ is semiprime, pick any $\mathbf{f} \in M$ such that 
	$f_i A \mathbf{f} \subseteq N$ for all $i$.
	It follows that $(f_i^\ast \mathbf{f})_j^\ast f_i^\ast \mathbf{f} =f_j^\ast f_i f_i^\ast \mathbf{f} \in A f_i A \mathbf{f} \subseteq N$ for all $i,j$. Since $N$ is real, it follows that $f_i^\ast \mathbf{f} \in N$ for every $i$. By using again that $N$ is real, it follows that $\mathbf{f} \in N$.
	
	To prove that the ideal $(N : M)$ is real, pick any $a_1,\ldots,a_s \in A$ such that
	$\sum_{i=1}^s a_i^\ast a_i \in (N : M)$. 
	For every $\mathbf{m} \in M$ and every $j$ we have
	$$\sum_{i=1}^s (a_i \mathbf{m})_j^\ast a_i \mathbf{m}= 
	m_j^\ast \sum_{i=1}^s  a_i^\ast a_i \mathbf{m} \in N$$
	Since $N$ is real, it follows that $a_i \mathbf{m}\in N$ for every $i=1,\ldots,s$
	and every $\mathbf{m} \in M$. Thus $a_1,\ldots,a_s \in (N : M)$.
	
	$(4) \Rightarrow (5)$ Suppose that $\mathfrak{p}$ is a two-sided real prime ideal.
	To prove that the submodule $C_{\mathfrak{p},\mathbf{u}}$ is real,
	pick $\mathbf{f}_1,\ldots,\mathbf{f}_m \in A^n$ such that
	$$\sum_{i=1}^m f_{i,j}^\ast \mathbf{f}_i \in C_{\mathfrak{p},\mathbf{u}} \text{ for every } j=1,\ldots,n.$$
	By the definition of $C_{\mathfrak{p},\mathbf{u}}$, it follows that 
	$$\langle \sum_{i=1}^m f_{i,j}^\ast \mathbf{f}_i,\mathbf{u} \rangle \in \mathfrak{p} \text{ for every } j=1,\ldots,n.$$
	If we multiply each equation by $u_j^\ast$ and sum over $j$ we get 
	$$\sum_{i=1}^m \langle \mathbf{f}_i,\mathbf{u} \rangle^\ast \langle \mathbf{f}_i,\mathbf{u} \rangle=\sum_{j=1}^n u_j^\ast 
	\sum_{i=1}^m f_{i,j}^\ast \langle \mathbf{f}_i,\mathbf{u} \rangle 
	=\sum_{j=1}^n u_j^\ast \langle \sum_{i=1}^m f_{i,j}^\ast \mathbf{f}_i,
	\mathbf{u} \rangle \in \mathfrak{p}$$
	Since $\mathfrak{p}$ is real, it follows that  
	$\langle \mathbf{f}_i,\mathbf{u} \rangle \in \mathfrak{p}$ for every $i$. 
	Therefore $\mathbf{f}_i \in C_{\mathfrak{p},\mathbf{u}}$ for every $i$. 
	Since $\mathfrak{p}$ is two-sided and prime, we have by Example \ref{cpu}
	that $C_{\mathfrak{p},\mathbf{u}}$ is also a $\mathfrak{p}$-prime submodule.
	
	$(5) \Rightarrow (1)$ follows from the observation that the intersection of a family of real 
	submodules is also a real submodule.
	
	For the rest of the proof we assume that $A$ is commutative.
	
	$(3) \Rightarrow (4)$ follows from the formula \eqref{kcf} in claim (4) of Proposition \ref{class}.
	
	$(2) \Rightarrow (3)$ Let $N$ be a semiprime submodule of $M$ 
	such that the ideal $(N : M)$ is real. By Theorem \ref{intthm1}, 
	$$N=\bigcap_{\mathfrak{p} \in \mathrm{Spec A}} K(N,\mathfrak{p}).$$ 
	By claim (4) of Proposition \ref{class}, we can omit those $\mathfrak{p}$ 
	that do not contain $(N : M)$ (because $K(N,\mathfrak{p})=M$ in this case).
	Clearly, we can also omit those $\mathfrak{p}$ that are not minimal over $(N : M)$.
	Since the ideal $(N : M)$ is real, we know by \cite[Lemma 2.9]{lam} that
	every minimal prime ideal over it is also real.
\end{proof}

\begin{rem}
	Let us show that for $R=\RR[x,y]$ there exists a submodule $N$ of $R^2$ such that 
	the ideal $(N : R^2)$ is real but $N$ is not real.
	Let $N$ be the submodule of $R^2$ generated by $x(x,y)$ and $y(x,y)$.
	Note that $N$ is not semiprime since it does not contain $(x,y)$. So $N$ is not real.
	A short computation shows that $(N : R^2)=\{0\}$ which is a real ideal of $R$.
	Namely, if $a \in (N : R^2)$ then $a(1,0) \in N$. By the definition of $N$ there
	exist $u,v \in R$ such that $a(1,0)=ux(x,y)+vy(x,y)$. It follows that 
	$a=(ux+vy)x$ and $0=(ux+vy)y$ which implies that $a=0$.
\end{rem}

As a corollary of Theorem \ref{realthm1},
we obtain a vector version of Theorem \ref{thma}.
In the following, $\mathbf{f},\mathbf{g}_i$ are row vectors 
and $\mathbf{u},\mathbf{v}$ are column vectors.

\begin{cor}
	\label{realthm2}
	Let $R=\RR[x_1,\ldots,x_d]$.
	For every $\mathbf{g}_1,\ldots,\mathbf{g}_m, \mathbf{f} \in R^n$ the following claims are equivalent:
	\begin{enumerate}
		\item[(1)] For every $a \in \RR^d$ and $\mathbf{v} \in \RR^n$ such that 
		$\mathbf{g}_1(a)\mathbf{v}=\ldots=\mathbf{g}_m(a)\mathbf{v}=0$
		we have $\mathbf{f}(a)\mathbf{v}=0$.
		\item[(2)] $\mathbf{f}$ belongs to the smallest real submodule of $R^n$ which contains $\mathbf{g}_1,\ldots,\mathbf{g}_m$.
	\end{enumerate}
\end{cor}

\begin{proof}
By Theorem \ref{thma}, every real ideal $\mathfrak{p}$ of $R$ satisfies
\begin{equation} \label{cor2eq1}
\mathfrak{p}=\bigcap_{\mathfrak{p} \subseteq \mathfrak{m}_a} \mathfrak{m}_a
\end{equation}
where $a$ runs through all real points and $\mathfrak{m}_a:=\{q \in R \mid q(a)=0\}$.
It follows that for every real prime ideal $\mathfrak{p}$ of $R$ and for every
$\mathbf{u} \in R^n \setminus \mathfrak{p}^n$
\begin{equation} \label{cor2eq2}
  	C_{\mathfrak{p},\mathbf{u}}
	=\bigcap_{\mathfrak{p} \subseteq \mathfrak{m}_a} C_{\mathfrak{m}_a,\mathbf{u}}
	=\bigcap_{\mathfrak{p} \subseteq \mathfrak{m}_a} C_{\mathfrak{m}_a,\mathbf{u}(a)}
\end{equation}
Let $N$ be the smallest real submodule of $R^n$ which contains $\mathbf{g}_1,\ldots,\mathbf{g}_m$.
By Theorem \ref{realthm1} and formula \eqref{cor2eq2}, it follows that
\begin{equation} \label{cor2eq3}
N=\bigcap_{N \subseteq C_{\mathfrak{p},\mathbf{u}}} C_{\mathfrak{p},\mathbf{u}}
=\bigcap_{N \subseteq C_{\mathfrak{m}_a,\mathbf{v}}} C_{\mathfrak{m}_a,\mathbf{v}} 
\end{equation}
where $a$ runs through $\RR^d$ and $\mathbf{v}$ runs through $\RR^n$. 
Therefore, any element $\mathbf{f} \in R^n$ belongs to $N$ iff
	every submodule of the  form  $C_{\mathfrak{m}_a,\mathbf{v}}$
	(with $a \in \RR^d,\mathbf{v} \in \RR^n$) which contains 
	$\mathbf{g}_1,\ldots,\mathbf{g}_m$ also contains $\mathbf{f}$.
\end{proof}

\section{Proof of Theorem \ref{thmd}}
\label{sec4}

The aim of this section is to deduce the matrix version of Theorem \ref{thma}
(i.e. Theorem \ref{thmd}) from the vector version of Theorem \ref{thma}
(i.e. from Corollary \ref{realthm2}.)

Recall that there is a natural 1-1 correspondence between the submodules of $A^n$ and the left ideals of $M_n(A)$. Recall also that a left ideal $I$ of a ring with involution $A$ is real
if for every $a_1,\ldots,a_r \in A$ such that $a_1 ^\ast a_1+\ldots+a_r^\ast a_r \in I$
we have $a_1,\ldots,a_r \in I$. 

\begin{prop}
	\label{realthm3}
	Let $A$ be an associative ring with involution.
	A submodule of $A^n$ is real iff the corresponding left ideal of $M_n(A)$ is real.
\end{prop}

\begin{proof}
	If $N$ is a submodule of $M:=A^n$ then the corresponding left ideal $I$ of $M_n(A)$
	consists of all matrices that have all rows in $N$.
	
	Suppose that $N$ is real.
	To prove that $I$ is real pick any  $A_1,\ldots,A_l \in M_n(A)$ 
	such that $A_1^\ast A_1+\ldots+A_l^\ast A_l \in I.$ Note that 
	$A_i^\ast A_i =\sum_{j=1}^n \mathbf{a}_{i,j}^\ast \mathbf{a}_{i,j}$ 
	where $\mathbf{a}_{i,1},\ldots,\mathbf{a}_{i,n}$ are the rows of $A_i$ 
	and so the $k$-th row of $A_i^\ast A_i$ is equal to 
	$\sum_{j=1}^n a_{i,j,k}^\ast \mathbf{a}_{i,j}$ where $a_{i,j,k}$ is the $k$-th component of $\mathbf{a}_{i,j}$. Since each row of $\sum_{i=1}^l A_i^\ast A_i$ 
	is in $N$, we have 
	$$\sum_{i=1}^l \sum_{j=1}^n a_{i,j,k}^\ast \mathbf{a}_{i,j} \in N$$ 
	for every $k$.
	The assumption that $N$ is real implies that $\mathbf{a}_{i,j} \in N$ for all $i,j$ 
	which implies that $A_i \in I$ for each $i$. Therefore $I$ is real.
	
	To prove the converse suppose that $I$ is a real.
	To prove that $N$ is real pick any $\mathbf{f}_1,\ldots,\mathbf{f}_k \in A^n$
	such that $\sum_{i=1}^n f_{i,j}^\ast \mathbf{f}_i \in N$ for every $j$.
	It follows that $\sum_{i=1}^n \mathbf{f}_i^\ast \mathbf{f}_i \in I$.
	Let $F_i$ be the matrix whose first row is $\mathbf{f}_i$ and other rows are zero.
	Since $\sum_{i=1}^n F_i^\ast F_i \in I$ and $I$ is real, it follows that $F_i \in I$
	for all $i$, and so $\mathbf{f}_i \in N$ for all $i$. Therefore $N$ is real.
\end{proof}

We are now ready for the proof of Theorem \ref{thmd}

\begin{proof}[Proof of Theorem \ref{thmd}]
	Write $A=\RR[x_1,\ldots,x_d]$ and pick $P_1,\ldots,P_m,Q \in M_n(A)$. 
	Let $\mathbf{p}_{1,1},\ldots,\mathbf{p}_{m,n}$ be the rows of $P_1,\ldots,P_m$ 
	and let $\mathbf{q}_1,\ldots,\mathbf{q}_n$ be the rows of $Q$.
	Let $N$ be the smallest real submodule of $A^n$ that contains all $\mathbf{p}_{i,j}$ 
	and let $I$ be the left ideal of $M_n(A)$ that corresponds to $N$. By Proposition 
	\ref{realthm3}, $I$ is the smallest real left ideal of $M_n(A)$ that contains all $P_i$.
	The following claims are equivalent:
	\begin{enumerate}
		\item $Q(a)\mathbf{v}=0$ for every $a \in \RR^d$ and every $\mathbf{v} \in \RR^n$ such that $P_i(a)\mathbf{v}=0$ for all $i$.
		\item $\mathbf{q}_1(a)\mathbf{v}=\ldots=\mathbf{q}_n(a)\mathbf{v}=\mathbf{0}$ 
		for every $a \in \RR^d$ and every $\mathbf{v} \in \RR^n$ such that 
		$\mathbf{p}_{i,j}(a)\mathbf{v}=0$ for all $i,j$.
		\item $N$ contains $\mathbf{q}_1,\ldots,\mathbf{q}_n$.
		\item $I$ contains $Q$.
	\end{enumerate}
	The equivalence of (2) and (3) follows from Corollary \ref{realthm2}.
\end{proof}

\section{Submodules of $(D_R)^n$}
\label{sec5}

The aim of this and the next section is to make technical preparations for the proof
of Theorem \ref{thme}. More precisely, we extend the results of Section \ref{sec2}
from commutative rings to a class of non-commutative associative rings that includes quaternionic polynomials. 

Let us first introduce the class of rings that we are going to study.
Let $D$ be a finite-dimensional central division algebra over a field $F$ and let $R$ be a commutative $F$-algebra. Write $D_R= D \otimes_F R$ and note that this is an $R$-algebra. 
For example, $\HH_{\RR[x_1,\ldots,x_d]}=\HH_c[x_1,\ldots,x_d]$.

By \cite[Lemma 17.4]{gwbook}, the mappings $\mathfrak{a} \mapsto \mathfrak{a}D_R=D \otimes_R \mathfrak{a}$ and $I \mapsto I \cap R$ give a 1-1 correspondence between the ideals of $R$ 
and the two-sided ideals of $D_R$.
Moreover,  prime ideals of $R$ correspond to two-sided prime ideals of $D_R$.
By the proof of \cite[Lemma 9.4]{pierce}, every $F$-basis of $D$ is also an $R$-basis of $D_R$.

\begin{lem}
	\label{dad}
	If $v_1,\ldots,v_l$ is an $F$-basis of $D$ and $c_1,\ldots,c_l \in R$
	then every two-sided ideal of $D_R$ which contains $\sum_{i=1}^l c_i v_i$
	also contains $c_1,\ldots,c_l$.
\end{lem}

\begin{proof}
	For each $i=1,\ldots,l$ let $\phi_i$ be an $F$-linear endomorphism of $D$
	which sends $v_i$ to $1$ and all other $v_j$ to $0$. By \cite[Corollary 3]{qsatz1}
	there exist elements $u_{i,k} \in D$ such that 
	$\phi_i(x)=\sum_{k=1}^l u_{i,k} x v_k$ for every $x \in D$.
	So,
	\begin{equation*}
		\label{q1c3}
		c_i=\sum_{j=1}^l c_j \phi_i(v_j)=
		\sum_{k=1}^l u_{i,k} \left( \sum_{j=1}^l c_j v_j \right) v_k \in
		\sum D \left( \sum_{j=1}^l c_j v_j \right) D
	\end{equation*}
	which is contained in every two-sided ideal of $D_R$ that contains $\sum_{j=1}^l c_j v_j$.
\end{proof}

Pick $n \in \NN$. We can consider $(D_R)^n$ either as a (left) $D_R$-module or as an $R$-module. In the first case we denote it by $M$ and in the second case by $\tilde{M}$.
Let $\mathbf{e}_1,\ldots,\mathbf{e}_n$ be the standard basis of $M$ and let
$v_1,\ldots,v_l$ be an $F$-basis of $D$.
Then $v_1 \mathbf{e}_1,\ldots,v_l \mathbf{e}_1,\ldots,v_1 \mathbf{e}_n,\ldots,v_l \mathbf{e}_n$ is a basis of $\tilde{M}$. 
Every submodule of $M$ is also a submodule of $\tilde{M}$.

\begin{prop}
	\label{primme}
	Let $N$ be a submodule of $M$ and $\mathfrak{p}$ a prime ideal of $R$.
	\begin{enumerate}
		\item $N$ is a semiprime submodule of $M$ iff it is a semiprime submodule of $\tilde{M}$.
		\item $N$ is a prime submodule of $M$ iff it is a prime submodule of $\tilde{M}$.
		\item $N$ is $\mathfrak{p}D_R$-prime as a submodule of $M$
		iff it is $\mathfrak{p}$-prime as a submodule of $\tilde{M}$.
		\item The smallest $\mathfrak{p}D_R$-prime submodule of $M$ that contains $N$
		is equal to the smallest $\mathfrak{p}$-prime submodule of $\tilde{M}$ that contains $N$.
	\end{enumerate}
\end{prop}

\begin{proof}
	(1) Suppose that $N$ is a semiprime submodule of $M$. Pick any 
	$$\mathbf{f}=\sum_{i=1}^n f_i \mathbf{e}_i=\sum_{i=1}^n \sum_{j=1}^l c_{i,j} v_j \mathbf{e}_i \in M$$
	such that $c_{i,j} \mathbf{f} \in N$ for every $i$ and $j$. It follows that
	$$f_i a \mathbf{f}=\sum_{j=1}^l c_{i,j} v_j a \mathbf{f} =\sum_{j=1}^l v_j a c_{i,j} \mathbf{f}\in N$$
	for every $a \in D_R$ and every $i$. Since $N$ is semiprime in $M$, it follows that $\mathbf{f} \in N$. So, $N$ is semiprime in $\tilde{M}$.
	
	Suppose that $N$ is a semiprime submodule of $\tilde{M}$.
	To prove that $N$ is semiprime in $M$
	pick any $\mathbf{f}=\sum_{i=1}^n f_i \mathbf{e}_i \in M$ 
	such that $f_i  D_R \mathbf{f} \subseteq N$ for every $i$.
	Write $f_i=\sum_{j=1}^l c_{i,j}v_j$ for every $i$ and recall that
	$c_{i,j} \in \sum D_R f_i D_R$ for every $i$ and $j$ by Lemma \ref{dad}.
	It follows that $c_{i,j} \mathbf{f} \in \sum D_R f_i D_R \mathbf{f} \subseteq N$ 
	for every $i$ and $j$.
	Since $N$ is semiprime in $\tilde{M}$, it follows that $\mathbf{f} \in N$.
	
	(2) Suppose that $N$ is a prime submodule of $M$. Pick $r \in R$ and $\mathbf{m} \in M$ such that $r \mathbf{m} \in N$. It follows that $r D_R \mathbf{m} \subseteq N$. 
	Since $N$ is prime, it follows that either $rM \subseteq N$ or $\mathbf{m} \in N$. Therefore $N$ is a prime submodule of $\tilde{M}$.
	
	Suppose that $N$ is a prime submodule of $\tilde{M}$. To show that $N$
	is a prime submodule of $M$, pick $a \in D_R$ and $\mathbf{m} \in M$ such that
	$a D_R \mathbf{m} \subseteq N$. Write $a=\sum_{j=1}^l \alpha_j v_j$ and recall that
	$\alpha_j \in \sum D_R a D_R$ for $j=1,\ldots,l$ by Lemma \ref{dad}.
	Therefore $\alpha_j \mathbf{m} \in \sum D_R a D_R \mathbf{m} \subseteq N$ for $j=1,\ldots,l$. 
	Since $N$ is a prime submodule of $\tilde{M}$, it follows that either $\mathbf{m} \in N$
	or $\alpha_j M \subseteq N$ for $j=1,\ldots,l$. The latter implies that $a M \subseteq N$.
	Thus $N$ is a prime submodule of $M$.
	
	(3) By the definition of $\tilde{M}$ we have $(N : \tilde{M})=(N : M) \cap R$.
	By the correspondence between the ideals of $R$ and the two-sided ideals
	of $D_R$, it follows that $(N : \tilde{M})D_R=(N : M)$. In particular,
	$(N : \tilde{M})=\mathfrak{p}$ iff $(N : M)=\mathfrak{p} D_R$.
	
	(4) Let $K_{\tilde{M}} (N,\mathfrak{p})$ be the smallest $\mathfrak{p}$-prime submodule of 
	$\tilde{M}$ that contains $N$. Since the $R$-linear isomorphism of $(D_R)^n$ to $R^{nl}$ 
	sends $\mathfrak{p}(D_R)^n=(\mathfrak{p}D_R)^n$ to $\mathfrak{p}^{nl}$ we have by 
	claim (3) of Proposition \ref{class} that
	\begin{equation}
		\label{ktildem}
		K_{\tilde{M}} (N,\mathfrak{p})=\{\mathbf{m} \in (D_R)^n \mid 
		\exists c \in R\setminus \mathfrak{p} \colon c \mathfrak{m} \in N+\mathfrak{p}(D_R)^n\}.
	\end{equation}
	To show that $K_{\tilde{M}}(N,\mathfrak{p})$ is a submodule of $M$ 
	pick any $\mathbf{m} \in K_{\tilde{M}}(N,\mathfrak{p})$ and any $a \in D_R$.
	By formula \eqref{ktildem}, there exists $c \in R \setminus \mathfrak{p}$ 
	such that $c \mathbf{m} \in N+\mathfrak{p}(D_R)^n$.
	Since  $N$ and $\mathfrak{p}(D_R)^n$ are submodules of $M$,
	it follows that $c a \mathbf{m} \in N+\mathfrak{p}(D_R)^n$ which implies that 
	$a \mathbf{m} \in K_{\tilde{M}}(N,\mathfrak{p})$.
	By claim (3), 
	$K_{\tilde{M}} (N,\mathfrak{p})$ is a $\mathfrak{p}D_R$-prime submodule of $M$.
	By its definition, it is also the smallest 
	$\mathfrak{p}D_R$-prime submodule of $M$ that contains $N$.
\end{proof}

We can now extend the implication $(1) \Rightarrow (2)$ of Theorem \ref{intthm1}
from commutative rings to rings of the form $A=D_R$.

\begin{thm}
	\label{ncintthm2}
	Every semiprime submodule $N$ of $M=(D_R)^n$ satisfies
	$$N=\bigcap_P K(N,P)$$
	where $P$ runs through all two-sided prime ideals of $D_R$.
	We can omit $P$ from the intersection if $K(N,P)=M$
	which holds iff $P \not\supseteq (N : M)$.
\end{thm}

\begin{proof}
	Since $N$ is a semiprime submodule of $M$, it is also a semiprime submodule of
	$\tilde{M}$ by claim (2) of Proposition \ref{primme}.
	By the implication $(1) \Rightarrow (2)$ of Theorem \ref{intthm1}, we have
	$N=\bigcap_{\mathfrak{p}\in \mathrm{Spec}R} K_{\tilde{M}}(N,\mathfrak{p})$.
	By claim (4) of Proposition \ref{primme}, $K_{\tilde{M}}(N,\mathfrak{p})=K(N,\mathfrak{p}D_R)$.
	
	By claim (4) of Proposition \ref{class},  $K_{\tilde{M}}(N,\mathfrak{p})=\tilde{M}$
	iff $\mathfrak{p} \not\supseteq (N : \tilde{M})$. 
	The correspondence between the ideals of $R$ and the two-sided ideals of $D_R$ sends
	$\mathfrak{p}$ to $\mathfrak{p}D_R$ and $(N : \tilde{M})=(N : M) \cap R$ to $(N : M)$
	preserving inclusions. It follows that $K(N,\mathfrak{p}D_R)=M$ iff 
	$\mathfrak{p}D_R\not\supseteq (N : M)$.
\end{proof}

\begin{rem}
	\label{rem7}
	Note that formula \eqref{ktildem} implies that for every prime $P$
	\begin{equation*}
		K(N,P)=\{\mathbf{m} \in (D_R)^n \mid 
		\exists c \in D_R \setminus P \colon 
		c D_R \mathbf{m} \subseteq N+P^n\}
	\end{equation*}
	We do not know how to extend this formula from $D_R$ to general rings.
\end{rem}

\section{Maximal submodules of $(D_R)^n$.}
\label{sec6}

In this section we finish the proof of the promised non-commutative version
of Theorem \ref{intthm1}. We have already proved that the implication
$(1) \Rightarrow (2)$ of Theorem \ref{intthm1} extends from $R$ to $D_R$.
It remains to show that the implication $(2) \Rightarrow (3)$ 
also extends from $R$ to $D_R$.

Recall that Proposition \ref{class} gives a classification of all prime submodules of 
$R^n$ where $R$ is commutative.
Proposition \ref{ncclass} will extend this classification to all prime submodules of $(D_R)^n$ where $D$ is a finite-dimensional central division $F$-algebra, $R$
is a commutative $F$-algebra and $D_R=D \otimes_F R$.

\begin{prop} 
	\label{ncclass}
	Let $D$ and $R$ be as above and let $\mathfrak{p} \ne R$ be a prime ideal of $R$.
	There is a $1-1-1-1$ correspondence between
	\begin{enumerate}
		\item $\mathfrak{p}D_R$-prime submodules of $(D_R)^n$,
		\item $0$-prime submodules of $(D_{\overline{R}})^n$ where $\overline{R}=R/\mathfrak{p}$ and $D_{\overline{R}}=D \otimes_F \overline{R}$.
		\item proper submodules of $(D_k)^n$ where $k$ is the field of fractions of $\overline{R}$ and $D_k=D \otimes_F k$,
		\item proper subspaces of $\tilde{D}^{nr}$ where a natural number $r$ and a
		division algebra $\tilde{D}$ are such that $D_k \cong M_r(\tilde{D})$. 
		(Note that $D_k$ is simple by \cite[Lemma 17.4]{gwbook}
		and Artinian since $\dim_k D_k=\dim_F D$.)
	\end{enumerate}
\end{prop}

\begin{proof} 
	$(1) \leftrightarrow (2)$ Let $\pi \colon R \to \overline{R}$ be the canonical homomorphism. Write $\phi=\id_D \otimes \pi \colon D_R \to D_{\overline{R}}$
	and $\bphi = \phi^n \colon M \to \overline{M}$ where 
	$M=(D_R)^n$ and $\overline{M}=(D_{\overline{R}})^n=\bphi(M)$.
	We will show that $N \to \bphi(N)$ and $\overline{N} \to \bphi^{-1}(\overline{N})$
	give a $1-1$ correspondence between (1) and (2).
	
	Since $\bphi$ is onto, we have $\bphi(\bphi^{-1}(\overline{N}))=\overline{N}$ for every
	$\overline{N}$. On the other hand $\bphi^{-1}(\bphi(N))=N+\ker \bphi$ for every $N$
	and $\ker \bphi =(\mathfrak{p} D_R)^n \subseteq N$ by Example \ref{smallest}.
	
	Suppose that $N$ is $\mathfrak{p}D_R$-prime. To show that $\overline{N}:=\bphi(N)$ is $0$-prime,
	pick $\phi(x) \in (\overline{N} : \overline{M})$. It follows that $\phi(x) \bphi(M) \subseteq \bphi(N)$. Therefore $xM \subseteq N+\ker(\bphi)$. As above, $\ker \bphi \subseteq N$.
	Therefore $x \in (N:M)=\mathfrak{p} D_R$ and so $\phi(x)=0$.
	
	Suppose that $\overline{N}$ is $0$-prime. 
	To show that $N:=\bphi^{-1}(\overline{N})$ is $\mathfrak{p}D_R$-prime,
	pick $x \in D_R$ and $\mathbf{m} \in M$ such that $x D_R \mathbf{m} \subseteq N$.
	It follows that $\phi(x) \phi(D_R) \bphi(\mathbf{m}) \subseteq \overline{N}$. 
	Since $\overline{N}$ is $0$-prime, it follows that either 
	$\phi(x) \in (\overline{N} : \overline{M})=\{0\}$ or 
	$\bphi(\mathbf{m}) \in \overline{N}$, which implies that 
	either $x \in \mathfrak{p}D_R$ or $\mathbf{m} \in N$. 
	Note also that $\mathfrak{p}D_R=\phi^{-1}(0)=
	\phi^{-1}((\overline{N} : \overline{M}))=(N : M)$.
	
	$(2) \leftrightarrow (3)$ We will show that $\overline{N} \to (\overline{R} \setminus \{0\})^{-1} \overline{N}$ and $N \to N \cap \overline{M}$ give a $1-1$ correspondence between (2) and (3).
	
	If $\overline{N}$ is a submodule of $\overline{M}$ then clearly 
	$(\overline{R} \setminus \{0\})^{-1} \overline{N}=k \overline{N}=D_k \overline{N}$ 
	is a submodule of $(D_k)^n$. Moreover, if $\overline{N}$ is $0$-prime then by
	clearing denominators we see that $(\overline{R} \setminus \{0\})^{-1} \overline{N}$
	is a proper submodule of $(D_k)^n$.
	
	Let $N$ be a proper submodule of $(D_k)^n$. Then $(N : (D_k)^n)$ is
	a proper two-sided ideal of $D_k$ and so it is trivial since $D_k$ is a simple ring 
	by \cite[Lemma 17.4]{gwbook}. Write $\overline{N}:=N \cap \overline{M}$ and
	note that $(\overline{N} : \overline{M})$ is a subset of $(N : (D_k)^n)=\{0\}$.
	To show that $\overline{N}$ is prime pick $z \in \overline{R}$ and 
	$\mathbf{m} \in \overline{M}$ such that 
	$z D_{\overline{R}} \mathbf{m}\subseteq \overline{N}$.
	If $z \ne 0$ then $\mathbf{m}\in z^{-1}\overline{N} \subseteq N$ which implies that
	$\mathbf{m} \in N \cap \overline{M}=\overline{N}$.
	
	To show that the correspondence is 1-1, note that $(\overline{R} \setminus \{0\})^{-1}(N \cap \overline{M})$ is clearly a subset of $N$ and that the opposite inclusion follows by
	putting all terms of all components of any element of $N$ to a common denominator. Note also that $\overline{N} \subseteq ((\overline{R} \setminus \{0\})^{-1} \overline{N}) \cap \overline{M}$ and that the opposite inclusion follows from the assumption that $\overline{N}$ is $0$-prime.
	
	$(3) \leftrightarrow (4)$ Let $N$ be a $ D_k$-submodule of $(D_k)^n$, 
	i.e. a $M_r(\tilde{D})$ submodule of $M_r(\tilde{D})^n=M_{r, rn}(\tilde{D})$. 
	Let $\Phi(N)$ be the set of all rows of all elements of $N$.
	Clearly, $\Phi(N)$ is a $\tilde{D}$-submodule of $\tilde{D}^{rn}$.
	Conversely, let $\tilde{N}$ be a $\tilde{D}$-submodule of $\tilde{D}^{rn}$.
	Let $\Psi(\tilde{N})$ be the set of all elements of $M_{r, rn}(\tilde{D})$ that have all rows in $\tilde{N}$. Clearly, $\Psi(\tilde{N})$ is a $M_r(\tilde{D})$ submodule
	of $M_{r, rn}(\tilde{D})$. A short computation shows that $\Phi$ and $\Psi$ are inverse to each other.
\end{proof}

\begin{rem}
\label{nclinalg}
Let $K$ be a skew-field and $m$ a natural number. There is a correspondence
between left and right subspaces of $K^m$ given by
\begin{gather*}
S \mapsto S^\circ:=\{\mathbf{m} \in K^m \mid \langle \mathbf{s},\mathbf{m} \rangle=0
\text{ for all } \mathbf{s} \in S\}\\
T \mapsto T_\circ:=\{\mathbf{m} \in K^m \mid \langle \mathbf{m},\mathbf{t} \rangle=0
\text{ for all } \mathbf{t} \in T\}
\end{gather*}
By the usual Gauss algorithm one can show that $(S^\circ)_\circ=S$ for 
every left subspace $S$ and $(T_\circ)^\circ=T$ for every right subspace $T$.
This correspondence is therefore 1-1 and clearly it reverses inclusions.

Maximal proper left subspaces correspond to minimal nontrivial
right subspaces. In other words, maximal proper left subspaces of $K^m$
coincide with subspaces of the form $C_{0,\mathbf{w}}=(\mathbf{w} K)_\circ$ where
$\mathbf{w} \in K^m \setminus\{\mathbf{0}\}$.

Furthermore, every proper left subspace $U$ of $K^m$ satisfies
$U=\bigcap_{\mathbf{w} \in U^\circ \setminus \{0\}} C_{0,\mathbf{w}}$.
If $U^\circ=\mathbf{w}_1 K+\ldots+\mathbf{w}_k K$ then 
$U=C_{0,\mathbf{w}_1} \cap \ldots \cap C_{0,\mathbf{w}_k}$.
\end{rem}

\begin{ex}
\label{excorr}
Let $R$, $D$ and $D_R$ be as above. Let $\mathfrak{p}$ be a prime ideal of $R$
and let $\mathbf{u}=(u_1,\ldots,u_n)$ be an element of $(D_R)^n \setminus (\mathfrak{p} D_R)^n$. By Example \ref{cpu}, $C_{\mathfrak{p}D_R,\mathbf{u}}$ 
is a $\mathfrak{p}D_R$-prime submodule of $(D_R)^n$. 
Let us find out to which submodule of $\tilde{D}^{rn}$ it corresponds. 
We need two steps.
	
Firstly, the submodule $C_{\mathfrak{p}D_R,\mathbf{u}}$ of $(D_R)^n$
corresponds to the	submodule $C_{0,\bar{\mathbf{u}}}$ of $M_r(\tilde{D})^n$ 
where $\bar{\mathbf{u}}=(\bar{u}_1,\ldots,\bar{u}_n)$ is the image of 
$\mathbf{u}$ under the natural map $(D_R)^n \to 
(D_{R/\mathfrak{p}})^n \to (D_k)^n \cong M_r(\tilde{D})^n$.
	
For each $i=1,\ldots,r$ let  $\mathbf{w}_i$ be the element of $\tilde{D}^{rn}$
whose first $r$ components are the entries of the $i$-th column of $\bar{u}_1$, the next
$r$ components are the entries of the $i$-th column of $\bar{u}_2$, etc. 
Then the submodule $C_{0,\bar{\mathbf{u}}}$ of $M_r(\tilde{D})^n$ corresponds to the subspace
$C_{0,\mathbf{w}_1} \cap \ldots \cap C_{0,\mathbf{w}_r}$ of $\tilde{D}^{rn}$. 

It follows that every maximal  $\mathfrak{p}D_R$-prime submodule of $(D_R)^n$ 
is of the form $C_{\mathfrak{p}D_R,\mathbf{u}}$ but the converse is false in general.
If $D_k$ is a division algebra (e.g. when $D=\HH$ and $\mathfrak{p}$ 
is a real prime ideal of $R$; see \cite[Proposition  1.6]{pierce}) 
then we also have the converse.
\end{ex}

Theorem \ref{ncmax} extends from $R$ to $D_R$ the implication $(2) \Rightarrow (3)$
of Theorem \ref{intthm1} and the implication $(3) \Rightarrow (4)$ of Theorem \ref{realthm1}.

\begin{thm}
	\label{ncmax}
	For every submodule $N$ of $(D_R)^n$ and every two-sided prime ideal $P$ of $D_R$ 
	we have
	$$K(N,P)=\bigcap_{N \subseteq C_{P,\mathbf{u}}}
	C_{P,\mathbf{u}}$$
	where $\mathbf{u}$ runs through $(D_R)^n \setminus P^n$.
\end{thm}

\begin{proof}
	By Example \ref{cpu}, $C_{P,\mathbf{u}}$ is a $P$-prime  submodule of $(D_R)^n$ and, by Example \ref{knp}, $K(N,P)$ is the smallest $P$-prime submodule of $(D_R)^n$ that contains $N$. 
	It follows that $N \subseteq C_{P,\mathbf{u}}$ iff
	$K(N,P) \subseteq C_{P,\mathbf{u}}$.
	
	Let $k$ be the field of fractions of $R/(P \cap R)$ and let
	a number $r$ and a division algebra $\tilde{D}$ be such that
	$D \otimes_R k \cong M_r(\tilde{D})$. Let $\tilde{N}$ be the
	subspace of $\tilde{D}^{nr}$ that corresponds to $K(N,P)$.
	By Remark \ref{nclinalg}, $\tilde{N}$ is equal to the intersection of all
	subspaces of the form $C_{0,\mathbf{w}}$ that contain it.
	By Example \ref{excorr}, each subspace of the form $C_{0,\mathbf{w}}$
	corresponds to a submodule of $(D_R)^n$ of the form 
	$C_{P,\mathbf{u}}$ for some very special $\mathbf{u}$.
\end{proof}

\section{Proof of Theorem \ref{thme}}
\label{sec7}

The first aim of this section is to complete the proof of the promised non-commutative version
of Theorem \ref{realthm1}; see Theorem \ref{qthm1}. As a corollary we obtain a vector
version of Theorem \ref{thmb}; see Corollary \ref{qthm2}. Finally, we deduce
Theorem \ref{thme} from Corollary \ref{qthm2}.

Let $D$ be a finite-dimensional central division algebra over a field $F$. Let $\ast$ be
an $F$-linear involution on $D$ (i.e. an involution of the first kind). We can extend this involution to $D_R$ so that it is trivial on $R$. 

It is well-known that in commutative rings 
a minimal prime ideal over a real ideal is always real; see \cite[Lemma 2.9]{lam}.
We will extend this observation from $R$ to $D_R$.

\begin{lem}
	\label{ncminprime}
	Let $I$ be an ideal of $R$ such that the ideal $D \otimes I$ is real.
	If $\mathfrak{p}$ is a minimal prime ideal over $I$ then the ideal
	$D \otimes \mathfrak{p}$  is also real.
\end{lem}

\begin{proof}
	Clearly the factor ring $(D \otimes R)/(D \otimes I)$ is isomorphic to $D \otimes (R/I)$.
	Therefore, we can assume that $I=0$ and that $\mathfrak{p}$ is a minimal prime ideal
	of $R$. It remains to show that if the ring $D \otimes R$ is real then the ring
	$(D \otimes R)/(D \otimes \mathfrak{p}) \cong D \otimes (R/\mathfrak{p})$ is real.
	(A ring with involution is \textit{real} if  
	$\sum_i a_i^\ast a_i=0$ implies that $a_i=0$ for all $i$.)
	
	\smallskip
	
	\textit{First step:} The ring $D \otimes R_{\mathfrak{p}}$ is real.
	
	Let $v_1,\ldots,v_l$ be a basis of $D$ and let 
	$a_i=\frac{r_{i,1}}{s}v_1+\ldots+\frac{r_{i,l}}{s}v_l \in D \otimes R_{\mathfrak{p}}$ be such that $\sum_i a_i^\ast a_i=0$. Write $b_i=r_{i,1} e_1+\ldots+r_{i,l}e_l$. 
	From $0=\sum_i a_i^\ast a_i=\frac{1}{s^2}\sum_i b_i^\ast b_i$ and $D \otimes R_{\mathfrak{p}} D \cong (D \otimes R)_{\mathfrak{p}}$, 
	it follows that there exists $t \in R\setminus \mathfrak{p}$ such that 
	$t \cdot \sum_i b_i^\ast b_i=0$. Thus $\sum_i (t  b_i)^\ast (t b_i)=0$. 
	Since the ring $D \otimes R$ is real, it follows that $t b_i=0$ for all $i$.
	It follows that $a_i=\frac{b_i}{s}=0$ for all $i$.
	
	\smallskip
	
	\textit{Second step:} $\mathfrak{p} R_{\mathfrak{p}}=0$.
	
	Recall that $\mathfrak{p} R_{\mathfrak{p}}$ is the only maximal ideal of $R_\mathfrak{p}$ and that there is a 1-1 correspondence between the prime ideals of $R_{\mathfrak{p}}$ and prime ideals of $R$ that are contained in $\mathfrak{p}$.
	Since $\mathfrak{p}$ is a minimal prime ideal of $R$, it follows that $\mathfrak{p} R_{\mathfrak{p}}$ is the only prime ideal of $R_{\mathfrak{p}}$. Since 
	$D \otimes R_{\mathfrak{p}}$ is real, it follows that $R_{\mathfrak{p}}$ is real.
	So $R_{\mathfrak{p}}$ is reduced, which implies that the intersection of all its prime ideals is zero. Therefore, $\mathfrak{p} R_{\mathfrak{p}}=0$.
	
	\smallskip
	
	\textit{Conclusion:} Let $k$ be the field of fractions of $R/\mathfrak{p}$. Clearly,
	$k=R_{\mathfrak{p}}/\mathfrak{p} R_{\mathfrak{p}}=R_\mathfrak{p}$ by the second step. Since $D \otimes k =D \otimes R_{\mathfrak{p}}$ is real by the first step,
	it follows that $D \otimes (R/\mathfrak{p})$ is real.
\end{proof}

\begin{rem}
\label{preserve}
For $\HH_c[x_1,\ldots,x_d]$, Lemma 2 is much easier to prove because 
the correspondence between the ideals of $R=\RR[x_1,\ldots,x_d]$ and
the two-sided ideals of $D_R=\HH_c[x_1,\ldots,x_d]$ preserves reality. 

Suppose that $J$ is a real ideal of $R$ and pick 
$w_s=\alpha_s+\beta_s i+\gamma_s j+\delta_s k \in D_R=R+R i+R j+R k$
such that $\sum_s w_s^\ast w_s \in J \cdot D_R=J+J i+J j+J k$.
It follows that $\sum_s (\alpha_s^2+\beta_s^2+\gamma_s^2+\delta_s^2) \in J$
which implies that $\alpha_s,\beta_s,\gamma_s,\delta_s \in J$ for all $s$.
Therefore, $w_s \in J \cdot D_R$ for all $s$. Conversely, if $I$ is a
two-sided real ideal of $D_R$ then clearly $I \cap R$ is a real ideal of $R$.
\end{rem}

We have already proved that the implication $(3) \Rightarrow (4)$ of Theorem
\ref{realthm1} holds for all rings of the form $D_R$; see Theorem \ref{ncmax}.
We  can now prove that the implication $(2) \Rightarrow (3)$ of Theorem \ref{realthm1}
also holds for all rings of the form $D_R$; see Theorem \ref{qthm1}.
Therefore, for $A=D_R$ with $D$ and $R$ as above, all five claims of Theorem \ref{realthm1} are equivalent.

\begin{thm}
\label{qthm1}
Let $D$ and $R$ be as above and let $N$ be a semiprime submodule of $M:=(D_R)^n$
such that  the ideal $(N : M)$ is real. Then
$$N=\bigcap_P K(N,P)$$ 
where $P$ runs through all two-sided real prime ideals of $D_R$.
\end{thm}

\begin{proof}
By Theorem \ref{ncintthm2}, $N=\bigcap_P K(N, P)$ where 
	$P$ runs through all two-sided prime ideals of $D_R$ that contain $(N : M)$. 
	We can omit those $P$ that are not minimal over $(N : M)$.
	By Lemma \ref{ncminprime}, the remaining $P$ are real. 
\end{proof}

We can now prove a vector version of Theorem \ref{thmb}. As in Corollary \ref{realthm2},
$\mathbf{f},\mathbf{g}_i$ are row vectors and $\mathbf{u},\mathbf{v}$ are column vectors.

\begin{cor}
	\label{qthm2}
	Let $A=\HH_c[x_1,\ldots,x_d]$.
	For every $\mathbf{g}_1,\ldots,\mathbf{g}_m, \mathbf{f} \in A^n$ the following are equivalent:
	\begin{enumerate}
		\item[(1)] For every $a \in \RR^d$ and $\mathbf{v} \in \HH^n$ such that 
		$\mathbf{g}_1(a)\mathbf{v}=\ldots=\mathbf{g}_m(a)\mathbf{v}=0$ we have $\mathbf{f}(a)\mathbf{v}=0$.
		\item[(2)] $\mathbf{f}$ belongs to the smallest real submodule of $A^n$ which contains $\mathbf{g}_1,\ldots,\mathbf{g}_m$.
	\end{enumerate}
\end{cor}

\begin{proof}
	By Theorem \ref{thmb}, every real prime ideal $P$ of $A$ satisfies
	\begin{equation} \label{cor3eq1}
	P=\bigcap_{P \subseteq \mathfrak{m}_a} \mathfrak{m}_a
	\end{equation}
	where $\mathfrak{m}_a=\{q \in A \mid q(a)=0\}$ for every $a \in \RR^d$. 
	From \eqref{cor3eq1} we obtain
	\begin{equation} \label{cor3eq2}
	C_{P,\mathbf{u}}=C_{\bigcap \mathfrak{m}_a,\mathbf{u}}
	=\bigcap_{P \subseteq \mathfrak{m}_a} C_{\mathfrak{m}_a,\mathbf{u}}
	=\bigcap_{P \subseteq \mathfrak{m}_a} C_{\mathfrak{m}_a,\mathbf{u}(a)}
    \end{equation}
    for every $P$ and $\mathbf{u}$. 	
    Let $N$ be the smallest real submodule of $A^n$ which contains $\mathbf{g}_1,\ldots,\mathbf{g}_m$.
     By Theorems  \ref{ncmax} and \ref{qthm1} and by equation \eqref{cor3eq2} we have
    \begin{equation} \label{cor3eq4}
	N=\bigcap_{N \subseteq C_{P,\mathbf{u}}} C_{P,\mathbf{u}}=\bigcap_{N \subseteq C_{\mathfrak{m}_a,\mathbf{v}}} C_{\mathfrak{m}_a,\mathbf{v}}
	\end{equation}
	where $a$ runs through $\RR^d$ and $\mathbf{v}$ runs through $\HH^n$. 
	In other words, an element $\mathbf{f} \in A^n$ 
	belongs to $N$ iff
	every submodule of the  form  $C_{\mathfrak{m}_a,\mathbf{v}}$
	(with $a \in \RR^d,\mathbf{v} \in \HH^n$) which contains 
	$\mathbf{g}_1,\ldots,\mathbf{g}_m$ also contains $\mathbf{f}$.
\end{proof}

Finally we deduce Theorem \ref{thme} from Corollary \ref{qthm2}
and Proposition \ref{realthm3}.

\begin{proof}[Proof of Theorem \ref{thme}]
	Let $\mathbf{p}_{1,1},\ldots,\mathbf{p}_{m,n}$ be the rows of $P_1,\ldots,P_m$
	and let $\mathbf{q}_1,\ldots,\mathbf{q}_n$ be the rows of $Q$. 
	The following claims are equivalent:
	\begin{enumerate}
		\item For every $a \in \RR^d$ and every $\mathbf{v} \in \HH^n$ such that $P_1(a)\mathbf{v}=\ldots=P_m(a)\mathbf{v}=0$, we have $Q(a)\mathbf{v}=0$.
		\item For every $a \in \RR^d$ and every $\mathbf{v} \in \HH^n$ such that 
		$\mathbf{p}_{i,j}(a) \mathbf{v}=0$ for all $i,j$ 
		we have $\mathbf{q}_k(a) \mathbf{v}=0$ for all $k$.
		\item All $\mathbf{q}_k$ belong to the smallest real submodule of 
		$\HH_c[x_1,\ldots,x_d]^n$ containing all $\mathbf{p}_{i,j}$. 
		\item $Q$ belongs to the smallest real left ideal of $M_n(\HH_c[x_1,\ldots,x_d])$ that contains $P_1,\ldots,P_m$.
	\end{enumerate}
	Clearly, (1) and (2) are equivalent.
	By Corollary \ref{qthm2}, (2) and (3) are equivalent.
	By Proposition \ref{realthm3}, (3) and (4) are equivalent.
\end{proof}

\end{document}